\documentclass[a4paper,DIV=12,11pt]{scrartcl}

\usepackage[utf8]{inputenc}
\usepackage{amsfonts,amssymb,amsmath}
\usepackage{epsfig}
\usepackage{MnSymbol}

\usepackage{graphics}
\usepackage{subfig}
\usepackage{color}
\usepackage[ruled,vlined,linesnumbered,noresetcount]{algorithm2e}

\usepackage{tikz}
\usetikzlibrary{plotmarks}
\usetikzlibrary{shapes,arrows}
\usetikzlibrary{positioning}
\usetikzlibrary{trees,mindmap}
\usetikzlibrary{fadings,shapes.arrows,shadows}
\usetikzlibrary{decorations.pathreplacing}
\usetikzlibrary{calc}

\renewcommand{\vec}{\boldsymbol}

\renewcommand{\llangle}{\left\langle}
\renewcommand{\rrangle}{\right\rangle}

\newcommand{\divergence}{\boldsymbol \nabla \cdot}
\newcommand{\gradient}{\boldsymbol \nabla}
\DeclareMathOperator*{\argmin}{\mathrm{arg\,min}}

\newcommand{\strain}{\boldsymbol \varepsilon}
\newcommand{\displacement}{\boldsymbol u}
\newcommand{\flux}{\boldsymbol q}

\newcommand{\permeability}{\boldsymbol \kappa}
\newcommand{\biotcoefficient}{\boldsymbol \alpha}

\newcommand{\uRHS}{\boldsymbol f^n_{\boldsymbol \theta,\Delta t}}
\newcommand{\pRHS}{h^n_{\boldsymbol \theta,\Delta t}}
\newcommand{\qRHS}{\boldsymbol g^n_{\boldsymbol \theta,\Delta t}}

\newcommand{\displacementSpace}{{\boldsymbol{V}}^n}
\newcommand{\pressureSpace}{{Q}^n}
\newcommand{\fluxSpace}{{\boldsymbol {W}}^n}
\newcommand{\displacementTestSpace}{{\boldsymbol{V}}_0}
\newcommand{\pressureTestSpace}{{Q}_0}
\newcommand{\fluxTestSpace}{{\boldsymbol{W}}_0}

\newcommand{\testdisplacement}{\boldsymbol v}
\newcommand{\testpressure}{q}
\newcommand{\testflux}{\boldsymbol w}

\newtheorem{defi}{Definition}[section]
\newtheorem{thm}[defi]{Theorem}
\newtheorem{rem}[defi]{Remark}
\newtheorem{cor}[defi]{Corollary}
\newtheorem{prob}[defi]{Problem}

\newenvironment{proof}{\paragraph{Proof.}}{\hfill$\square$}
\numberwithin{equation}{section}
\numberwithin{algocf}{section}

\begin{document}

\title{Iterative Coupling for Fully Dynamic Poroelasticity}

\author{M.\ Bause\thanks{bause@hsu-hh.de (corresponding author),
$^\dag$Jakub.Both@uib.no, $^\ddag$Florin.Radu@uib.no}, 
J.\ W.\ Both$^\dag$,  F.\ A.\ Radu$^\ddag$\\
{\small $^\ast$ Helmut Schmidt University, Faculty of 
Mechanical Engineering, Holstenhofweg 85,}\\ 
{\small 22043 Hamburg, Germany}\\
{\small ${}^\dag,{}^\ddag$ University of Bergen, Department of Mathematics, 
All\'{e}gaten 41},\\{\small 50520 Bergen, Norway}
}

\date{}

\maketitle

\begin{abstract}
\textbf{Abstract.} 
We present an iterative coupling scheme for the numerical approximation 
of the mixed hyperbolic-parabolic system of fully dynamic poroelasticity. We prove its 
convergence in the Banach space setting for an abstract semi-discretization in time that 
allows the application of the family of diagonally implicit Runge--Kutta methods. 
Recasting the semi-discrete solution as the minimizer of a properly defined energy 
functional, the proof of convergence uses its alternating minimization. The scheme is 
closely related to the undrained split for the quasi-static Biot system.
\end{abstract}



\section{Introduction}
\label{sec:introduction}

Information on flow in deformable porous media has become of increasing importance in 
various fields of natural sciences and technology. It offers an abundance of technical, 
geophysical, environmental and biomedical applications including modern material science 
polymers and metal foams, gaining significance particularly in lightweight design and 
aircraft industry, design of batteries or hydrogen fuel cells for green technologies,  
geothermal energy exploration or reservoir engineering as well as mechanism in the human 
body and food technology. Consequently, quantitative methods, based on numerical 
simulations, are desirable in analyzing experimental data and  designing theories based 
on 
mathematical concepts. Recently, the quasi-static Biot system (cf., e.g., 
\cite{MW13_BBR,S00_BBR}) has attracted researcher's interest and has been studied  as a 
proper model for the numerical simulation of flow in deformable porous media. The design, 
analysis and optimization of approximation techniques that are based on an iterative 
coupling of the subproblems of fluid flow and mechanical deformation were focused 
strongly. Iterative coupling offers the appreciable advantage over the fully coupled 
method that existing and highly developed discretizations and algebraic solver 
technologies can be reused. For the quasi-static Biot system, pioneering work is done in 
\cite{K11_BBR,MW13_BBR}. Further research is presented in, e.g., 
\cite{BRK17_2_BBR,BBNKR17_BBR,C15_BBR,C16_BBR,HKLP19_BBR,MWW14_BBR}.

In the case of larger contrast coefficients that stand for the ratio between 
the intrinsical characteristic time and the characteristic domain time scale the fully 
dynamic hyperbolic-parabolic system of poroelasticity has to be considered. In 
\cite{MW12_BBR}, this system (referred to as the Biot--Allard equations) is derived by 
asymptotic homogenization in the space and time variables. Here, to fix our ideas and 
carve out the key technique of proof, a simplified form of the system proposed in 
\cite{MW12_BBR} is studied. However, its mixed hyperbolic-parabolic structure is 
preserved.  Our modification of the fully dynamic poroelasticity model in  
\cite{MW12_BBR} 
comes through a simplication of the solution's convolution with the dynamic permeability  
that is defined as the spatial average of pore system Stokes solutions on the unit cell. 
The fully dynamic system of poroelasticity to be analyzed here is given by (cf.\ 
\cite[p.\ 
313]{S00_BBR})
\begin{subequations}
\label{Eq:BA_0_BBR}
\begin{align}
\label{Eq:BA_1_BBR}
\rho \, \partial_{t}^2 \boldsymbol{u}  - \divergence \left( \mathbb{C} 
\strain(\displacement) - \biotcoefficient p \right) &= \boldsymbol f\,,\\[1ex]
\label{Eq:BA_2_BBR}
\partial_t  \left(c_0 p + \biotcoefficient : \strain(\boldsymbol u)\right) + \divergence 
\flux &= h\,, \\[1ex]
\label{Eq:BA_3_BBR}
\permeability^{-1} \flux + \gradient p &= \boldsymbol g\,.
\end{align}
\end{subequations}
System \eqref{Eq:BA_0_BBR} is equiped with appropriate initial and boundary conditions. 
In 
 \eqref{Eq:BA_0_BBR}, the variable $\boldsymbol u$ is the unkown effective solid phase 
displacement and $p$ is the unkown effective pressure.  The quantity $\boldsymbol 
\varepsilon (\vec u)=(\nabla \vec u +(\nabla \vec u)^\top)/2$ denotes the symmetrized 
gradient or strain tensor. Further, $\rho$ is the effective mass density,  $\mathbb C$ is 
Gassmann's fourth order effective elasticity tensor, $\boldsymbol \alpha $ is Biot's 
pressure-storage coupling tensor and $c_0$ is the specific storage coefficient. In the 
three field formulation \eqref{Eq:BA_0_BBR}, the vector field $\boldsymbol q$ is Darcy's 
velocity and $\boldsymbol \kappa$ is the permeability tensor. All tensors are assumed to 
be symmetric, bounded and uniformly positive definite, the constants $\rho$ and $c_0$ are 
positive. By $\boldsymbol A : \boldsymbol B$ we denote the Frobenius inner product of 
$\boldsymbol A$ and $\boldsymbol B$. The functions on the right-hand side of 
\eqref{Eq:BA_0_BBR} are supposed to be elements in dual spaces and, therefore, can 
include 
body forces and surface data (boundary conditions).

So far, the numerical simulation of the system \eqref{Eq:BA_0_BBR} has been studied 
rarely in the literature despite its numerous applications in practice. This might be due 
to the mixed hyberbolic-parabolic character of the system and severe complexities 
involved 
in the construction of monolithic solver or iterative coupling schemes with guaranteed 
stability properties. Space-time finite element approximations of hyperbolic and 
parabolic 
problems and the quasi-static Biot system were recently proposed, analyzed and 
investigated numerically by the authors in \cite{BKRS18_BBR,BRK17_2_BBR,BRK17_BBR}. Here, 
we propose an iterative coupling scheme for the system \eqref{Eq:BA_0_BBR} and prove its 
convergence. This is done in Banach spaces for the semi-discretization in time of 
\eqref{Eq:BA_0_BBR}. An abstract setting is used for the time discretization such that 
the 
family of diagonally implicit Runge--Kutta methods becomes applicable.  The key 
ingredient of our proof of convergence is the observation that we can recast the 
semi-discrete approximation of \eqref{Eq:BA_0_BBR} as the minimizer of an energy 
functional in the displacement and Darcy velocity fields. To solve the minimization 
problem, the general and abstract framework of alternating minimization (cf.\ 
\cite{BKNR19_BBR,B19_BBR}) is applied. The resulting subproblems of this minimization are 
then reformulated as our final iterative coupling scheme. Thereby, the proof of 
convergence of the iterative scheme is traced back to the convergence of the alternating 
minimization approach. This shows that the latter provides an abstract and powerful tool 
of optimization for the design of iterative coupling schemes.

We use standard notation. In particular, we denote by $\langle\cdot, \cdot \rangle$ the 
standard inner product of $L^2(\Omega)$ and by $\|\cdot \|$ the norm of $L^2(\Omega)$.

\section{Variational formulation of a semi-discrete approximation of the system of 
dynamic poroelasticity}

Firstly, we discretize the continuous system of dynamic poroelasticity 
\eqref{Eq:BA_0_BBR} in time by using arbitrary (diagonally implicit) Runge--Kutta 
methods and formulate the semi-discrete approximation as solution to a minimization 
problem, following the approach in~\cite{BKNR19_BBR}. For this, we consider an 
equidistant partition $0=t_0 < t_1 < \ldots < t_N =T$ of the time interval of interest  
$[0,T]$ with time step size $\Delta t$. In the sequel, we use the following function 
spaces for displacement, pressure, and flux, respectively,
\begin{align*}
 \displacementSpace &:= 
 \left\{ \testdisplacement \in H^1(\Omega)^d \,\big| \, \testdisplacement^n \text{ 
satisfies prescribed BC at time }t_n \right\}\,,\\
 \pressureSpace &:= L^2(\Omega)\,,\\
 \fluxSpace &:= \left\{ \testflux \in H(\mathrm{div};\Omega) \, \big| \, \testflux\text{ 
satisfies prescribed BC at time }t_n \right\}\,.
\end{align*}
Further, let $\displacementTestSpace$, $\pressureTestSpace$, and $\fluxTestSpace$ denote 
the corresponding natural test spaces, and $\displacementTestSpace^\star$, 
$\pressureTestSpace^\star$, and $\fluxTestSpace^\star$ their dual spaces.

Applying any diagonally implicit Runge--Kutta method for the temporal discretization 
of~\eqref{Eq:BA_0_BBR}, eventually involves solving systems of the following structure.

\begin{prob}
In the $n$-th time step, find the displacement $\displacement^n \in  
\displacementSpace$, pressure $p^n\in \pressureSpace$, and flux $\flux^n \in\fluxSpace$, 
satisfying for all $(\testdisplacement,\testpressure, \testflux)\in 
\displacementTestSpace 
\times \pressureTestSpace \times \fluxTestSpace$ the equations
\begin{subequations}
\label{Eq:BA_semi-discrete-0_BBR}
\begin{align}
\label{Eq:BA_semi-discrete-1_BBR}
 \frac{\rho}{\Delta t^2} \llangle \displacement^n, \testdisplacement \rrangle + 
\theta_1 \llangle \mathbb{C} \strain(\displacement^n), \strain(\testdisplacement) 
\rrangle 
 - \llangle \biotcoefficient p^n, \strain(\testdisplacement) \rrangle &= \llangle \uRHS, 
\testdisplacement \rrangle\,,\\[1ex]
\label{Eq:BA_semi-discrete-2_BBR}
 c_0 \llangle p^n, \testpressure \rrangle + \llangle \biotcoefficient 
: \strain(\displacement^n), \testpressure \rrangle + \theta_2 \Delta t \llangle 
\divergence \flux^n, \testpressure \rrangle &= \llangle \pRHS, \testpressure 
\rrangle\,,\\[1ex]
\label{Eq:BA_semi-discrete-3_BBR}
 \llangle \permeability^{-1}\flux^n, \testflux \rrangle - \llangle p^n, 
\divergence \testflux \rrangle &= \llangle \qRHS, \testflux \rrangle\,.
\end{align}
\end{subequations}
\end{prob}

In \eqref{Eq:BA_semi-discrete-0_BBR},  the quantities $\theta_1,\theta_2\in (0,1]$ 
are discretization parameters, and the right-hand side functions  
$\uRHS\in\displacementTestSpace^\star$, $\pRHS \in \pressureTestSpace^\star$, $\qRHS \in 
\fluxTestSpace^\star$ include information on external volume and surface terms, as well 
as 
previous time steps depending on the choice of the implicit Runge--Kutta discretization. 

Assuming positive compressibility, i.e., $c_0>0$ for the specific storage coefficient, 
the semi-discrete approximation satisfies equivalently the following variational problem; 
cf.\ \cite{BKNR19_BBR} for the derivation of a similar equivalence in the framework of 
the quasi-static Biot system. 
\begin{prob}
Find $(\displacement^n, \flux^n) \in \displacementSpace \times \fluxSpace$, satisfying
\begin{align}
\label{Eq:BA_variational_BBR}
 (\displacement^n, \flux^n) &= \argmin_{(\displacement,\flux) \in \displacementSpace 
\times \fluxSpace} \mathcal{E}(\displacement,\flux)\,,
\end{align}
where the energy $\mathcal{E}:\displacementSpace \times \fluxSpace \rightarrow 
\mathbb{R}$ at time $t_n$ is defined by ($(\displacement,\flux)\in \displacementSpace 
\times \fluxSpace$)
\begin{equation}
\label{Def:EngFct}
\begin{aligned}
 \mathcal{E}(\displacement,\flux) &:= \frac{\rho}{2 \Delta t^2} \| \displacement \|^2 + 
\frac{\theta}{2} \llangle \mathbb{C} \strain(\displacement), \strain(\displacement) 
\rrangle 
 + \frac{\theta_1\theta_2\Delta t}{2} \llangle \permeability^{-1} \flux, \flux \rrangle 
\\[1ex] 
 &\qquad + \frac{\theta_1}{2c_0} \left\|\pRHS - \biotcoefficient : \strain(\displacement) 
- \theta_2 \Delta t \divergence \flux \right\|^2 - \llangle \uRHS, \displacement \rrangle 
- \llangle \qRHS, \flux \rrangle.
\end{aligned}
\end{equation}
\end{prob}

The semi-discrete pressure $p^n$ may then be recovered by the post-processing step
\begin{align}
\label{Eq:BA_pressure_BBR}
 p^n = c_0^{-1} \left( \pRHS - \biotcoefficient : \strain(\displacement^n) - \theta_2 
\Delta t \, \divergence \flux^n\right).
\end{align}

\section{Iterative coupling for the system of dynamic poroelasticity}
\label{sec:DefIterCoupl}

Following the philosophy of~\cite{BKNR19_BBR}, we propose an iterative coupling of 
the semi-discrete equations~\eqref{Eq:BA_semi-discrete-0_BBR} of dynamic poroelasticity  
by firstly applying the fundamental alternating minimization to the variational 
formulation~\eqref{Eq:BA_variational_BBR}; cf.\ Alg.~\ref{algorithm:undrained-split}. 
\begin{algorithm}[h]
 \caption{Single iteration of the alternating minimization}
 \label{algorithm:undrained-split}
 \normalsize
 \SetAlgoLined
 \DontPrintSemicolon
 
  \vspace{0.25em}
 
  Input: $(\displacement^{n,k-1},\flux^{n,k-1})\in\displacementSpace \times \fluxSpace$ \; 
\vspace{0.25em}
   
  Determine $\displacement^{n,k} := \argmin_{\displacement\in\displacementSpace}\, 
\mathcal{E}(\displacement,\flux^{n,k-1})$\; \vspace{0.25em}

  Determine $\flux^{n,k} := \argmin_{\flux\in\fluxSpace}\, 
\mathcal{E}(\displacement^{n,k},\flux)$\; \vspace{0.2em}
\end{algorithm}

Secondly, the resulting scheme is equivalently reformated in terms of a stabilized 
splitting scheme applied to the three-field 
formulation~\eqref{Eq:BA_semi-discrete-0_BBR}. 
For this, a pressure iterate $p^{n,k} = c_0^{-1} \left( \pRHS\right.$ $\left. - 
\biotcoefficient : 
\strain(\displacement^{n,k}) - \theta_2 \Delta t \, \divergence 
\flux^{n,k}\right)\in\pressureSpace$, $k\geq 0$, is introduced, consistent 
with~\eqref{Eq:BA_pressure_BBR}, and the optimality conditions corresponding to the two 
steps of Alg.~\ref{algorithm:undrained-split} are reformulated. The calculations are 
skipped here. We immediately present the resulting scheme, which 
in the end is closely related to the \textit{undrained split} for the quasi-static Biot 
system~\cite{K11_BBR}.

\begin{prob}
Let $(\displacement^{n,0},p^{n,0})\in\displacementSpace \times \pressureSpace$ be given 
and $k\geq 1$. 

\textbf{1.\ Step} (Update of mechanical deformation): For given 
$(\displacement^{n,k-1},p^{n,k-1})\in\displacementSpace \times \pressureSpace $, find 
$\displacement^{n,k}\in\displacementSpace$ satisfying for all $\testdisplacement \in 
\displacementTestSpace$,
\begin{align}
\label{Eq:BA_split_1_BBR}
  &\frac{\rho}{\Delta t^2} \llangle \displacement^{n,k}, \testdisplacement \rrangle + 
\theta_1 \llangle \mathbb{C} \strain(\displacement^{n,k}) + \frac{\biotcoefficient \otimes 
\biotcoefficient}{c_0} \strain(\displacement^{n,k} - \displacement^{n,k-1}), 
\strain(\testdisplacement) \rrangle\\
  \nonumber
  &\qquad- \theta_1 \llangle \biotcoefficient p^{n,k-1}, \strain(\testdisplacement) 
\rrangle = \llangle \uRHS, \testdisplacement \rrangle,
\end{align}
where $\otimes : \mathbb{R}^{d \times d} \times \mathbb{R}^{d \times d} \rightarrow 
\mathbb{R}^{d \times d \times d \times d}$ denotes the standard tensor product. 

\smallskip
\textbf{2.\ Step} (Update of Darcy velocity and pressure): For given 
$(\displacement^{n,k},p^{n,k-1})\in\displacementSpace \times \pressureSpace $ find 
$(p^{n,k},\flux^{n,k})\in \pressureSpace\times \fluxSpace$ satisfying for all 
$(\testpressure,\testflux) \in \pressureTestSpace \times \fluxTestSpace$, 
\begin{subequations}
\label{Eq:BA_split_0_BBR}
\begin{align}
\label{Eq:BA_split_2_BBR}
  c_0 \llangle p^{n,k}, \testpressure \rrangle  + \llangle \biotcoefficient 
: \strain(\displacement^{n,k}), \testpressure \rrangle  + \theta_2 \Delta t \llangle  
\divergence \flux^{n,k}, \testpressure \rrangle  &= \llangle \pRHS, \testpressure 
\rrangle\,,\\[1ex]
\label{Eq:BA_split_3_BBR}
 \llangle \permeability^{-1}\flux^{n,k}, \testflux \rrangle - \llangle p^{n,k}, 
\divergence \testflux \rrangle &= \llangle \qRHS, \testflux \rrangle\,.
\end{align}
\end{subequations}
\end{prob}

We note that the splitting scheme defined by~\eqref{Eq:BA_split_1_BBR}, 
\eqref{Eq:BA_split_0_BBR} utilizes the identical stabilization as the \textit{undrained 
split} for the quasi-static Biot equations~\cite{K11_BBR}.

\section{Convergence of the iterative coupling scheme}
\label{sec:ConvIterCoupl}

The identification of the undrained split approach~\eqref{Eq:BA_split_1_BBR}, 
\eqref{Eq:BA_split_0_BBR} as the application of the alternating minimization, cf.\ 
Alg.~\ref{algorithm:undrained-split}, to the variational 
problem~\eqref{Eq:BA_variational_BBR} yields the basis for a simple convergence 
analysis. For this, we utilize the following abstract convergence result, that is 
rewritten here in terms of the specific formulation of 
Alg.~\ref{algorithm:undrained-split}.

\begin{thm}[Convergence of the alternating 
minimization~\cite{B19_BBR}\label{Lemma:AM_convergence}]
Let $|\cdot|$, $|\cdot|_\mathrm{m}$, and $|\cdot|_\mathrm{f}$ denote semi-norms on 
$\displacementTestSpace \times \fluxTestSpace$, $\displacementTestSpace$, and 
$\fluxTestSpace$, respectively. Let $\beta_\mathrm{m},\beta_\mathrm{f}\geq 0$ satisfy the 
inequalities 
\begin{align*}
 |(\testdisplacement,\testflux)|^2 \geq \beta_\mathrm{m} | \testdisplacement 
|_\mathrm{m}^2\qquad \text{and} \qquad 
 |(\testdisplacement,\testflux)|^2 \geq \beta_\mathrm{f} | \testflux |_\mathrm{f}^2
\end{align*}
for all $(\testdisplacement,\testflux)\in\displacementTestSpace \times \fluxTestSpace$. 
Furthermore, assume that the energy functional $\mathcal{E}$ of \eqref{Def:EngFct} 
satisfies the following conditions:
\begin{itemize}
  \item The energy $\mathcal{E}$ is Frech\'et differentiable with $D\mathcal{E}$ denoting 
its derivative.

  \item The energy $\mathcal{E}$ is strongly convex wrt.\ $|\cdot|$ with modulus 
$\sigma>0$, i.e., for all $\displacement,\bar{\displacement}\in\displacementSpace$ and 
$\flux,\bar{\flux}\in\fluxSpace$ it holds that 
  \begin{align*}
  \mathcal{E}(\bar{\displacement},\bar{\flux}) \geq \mathcal{E}(\displacement,\flux) + 
\llangle D\mathcal{E}(\displacement,\flux), (\bar{\displacement} - 
\displacement,\bar{\flux} - \flux)  \rrangle + \frac{\sigma}{2} |(\bar{\displacement} - 
\displacement,\bar{\flux} - \flux)|^2\,.
  \end{align*}

  \item The partial functional derivatives $D_{\displacement} \mathcal{E}$ and 
$D_{\flux} \mathcal{E}$ are uniformly Lipschitz continuous wrt.\ 
$|\cdot|_\mathrm{m}$ and $|\cdot|_\mathrm{f}$ with Lipschitz constants $L_\mathrm{m}$ and 
$L_\mathrm{f}$, respectively, {i.e., for all $(\displacement,\flux) \in 
\displacementSpace 
\times \fluxSpace$ and $(\testdisplacement,\testflux) \in \displacementTestSpace \times 
\fluxTestSpace$ it holds that 
  \begin{align*}
   \mathcal{E}(\displacement + \testdisplacement, \flux) &\leq \mathcal{E}(\displacement, 
\flux) + \llangle D_{\displacement} \mathcal{E}(\displacement, \flux),\testdisplacement 
\rrangle + \frac{L_\mathrm{m}}{2} \left\| \testdisplacement 
\right\|_\mathrm{m}^2\,,\\[1ex]
   \mathcal{E}(\displacement, \flux + \testflux) &\leq \mathcal{E}(\displacement, \flux) + 
\llangle D_{\flux} \mathcal{E}(\displacement, \flux),\testflux \rrangle + 
\frac{L_\mathrm{f}}{2} \left\| \testflux \right\|_\mathrm{f}^2\,.
  \end{align*}}
\end{itemize}

\vspace*{-2ex} \noindent
Let $(\displacement^n,\flux^n)\in \displacementSpace \times \fluxSpace$ denote the 
solution to~\eqref{Eq:BA_variational_BBR}, and let $(\displacement^{n,k},\flux^{n,k})$ 
denote the corresponding approximation defined by Alg.~\ref{algorithm:undrained-split}. 
Then, for all $k\geq 1$ it follows that 
\begin{align}
\label{Eq:AM_convergence}
 &\mathcal{E}(\displacement^{n,k},\flux^{n,k}) - 
\mathcal{E}(\displacement^n,\flux^n)\\[1ex]
 \nonumber
 &\quad\leq \left( 1 - \frac{\beta_\mathrm{m} \sigma}{L_\mathrm{m}} \right) \left( 1 
- \frac{\beta_\mathrm{f} \sigma}{L_\mathrm{f}} \right) \, \left( 
\mathcal{E}(\displacement^{n,k-1},\flux^{n,k-1}) - \mathcal{E}(\displacement^n,\flux^n) 
\right)\,.
\end{align}
\end{thm}

A simple application of Theorem~\ref{Lemma:AM_convergence} now yields the main result of 
the work, namely the global linear convergence of the undrained 
split~\eqref{Eq:BA_split_1_BBR}, \eqref{Eq:BA_split_0_BBR}.

\begin{cor}[Linear convergence of the undrained 
split\label{Corollary:Convergence_UD}]
 Let $|\cdot|$ be defined by
 \begin{align*}
  \left| (\testdisplacement,\testflux) \right|^2 &:= \frac{\rho}{\Delta t^2} 
\left\| \testdisplacement \right\|^2 + \theta_1 \llangle \mathbb{C} 
\strain(\testdisplacement), \strain(\testdisplacement) \rrangle +\theta_1\theta_2 \Delta 
t 
\llangle \permeability^{-1} \testflux, \testflux \rrangle \\[0ex]
  &\qquad \qquad +\frac{\theta_1}{c_0} \left\| \biotcoefficient : 
\strain(\testdisplacement) + \theta_2 \Delta t \divergence \testflux \right\|^2
 \end{align*}
 for all $(\testdisplacement,\testflux) \in \displacementTestSpace \times \fluxTestSpace$. 
Furthermore, let $(\displacement^n,\flux^n)\in \displacementSpace \times \fluxSpace$ 
denote the solution to~\eqref{Eq:BA_variational_BBR}, and let 
$(\displacement^{n,k},\flux^{n,k})\in \displacementSpace \times \fluxSpace$ denote the 
corresponding approximation defined by Alg.~\ref{algorithm:undrained-split}. Then, for all 
$k\geq 1$ it holds that 
 \begin{align*}
  \left|(\displacement^{n,k} - \displacement^n, \flux^{n,k} - \flux^n)\right|^2 & \leq 
\left( \frac{\|\biotcoefficient :\mathbb{C}^{-1} : \biotcoefficient\|_\infty}{c_0 + 
\|\biotcoefficient :\mathbb{C}^{-1} : \biotcoefficient\|_\infty} \right)^2 \\[1ex]
  & \qquad  \cdot \left|(\displacement^{n,k-1} - \displacement^n, \flux^{n,k-1} - 
\flux^n)\right|^2.
 \end{align*}
\end{cor}

\begin{proof}
 We first examine convexity and smoothness properties of $\mathcal{E}$ defined in 
\eqref{Def:EngFct} by analyzing the second functional derivative of $\mathcal{E}$. For 
this, let $(\displacement,\flux)\in\displacementSpace\times \fluxSpace$ and 
$(\testdisplacement,\testflux)\in \displacementTestSpace \times \fluxTestSpace$ be 
arbitrary. Then, for the second functional derivative $D^2 
\mathcal{E}(\displacement,\flux):(\displacementTestSpace^\star \times 
\fluxTestSpace^\star)^2 \rightarrow \mathbb{R}$ of $\mathcal{E}$ it holds that
 \begin{align}
 \label{Eq:proof-aux:second-derivative-energy}
  \llangle D^2 \mathcal{E}(\displacement,\flux) (\testdisplacement,\testflux), 
(\testdisplacement,\testflux) \rrangle  = \left| (\testdisplacement,\testflux ) \right|^2.
 \end{align}
 Next, we define a norm $|\cdot|_\mathrm{m}$ on $\displacementTestSpace$ by considering 
the partial second functional derivative of $\mathcal{E}$ with respect to the displacement 
field,
 \begin{align*}
  \llangle D_{\displacement}^2 \mathcal{E}(\displacement,\flux) \testdisplacement, 
\testdisplacement \rrangle  
  =
  \frac{\rho}{\Delta t^2} \left\| \testdisplacement \right\|^2 + \theta_1 \llangle 
\mathbb{C} \strain(\testdisplacement), \strain(\testdisplacement) \rrangle 
+\frac{\theta_1}{c_0} \left\| \biotcoefficient : \strain(\testdisplacement)\right\|^2
  =: \left| \testdisplacement \right|_\mathrm{m}^2.
 \end{align*}
 Similarly, we define a norm $|\cdot|_\mathrm{f}$ on $\fluxTestSpace$ by means of 
 \begin{align*}
  \llangle D_{\flux}^2 \mathcal{E}(\displacement,\flux) \testflux, \testflux \rrangle  
  =
  \theta_1\theta_2 \Delta t \llangle \permeability^{-1} \testflux, \testflux \rrangle 
+\frac{\theta_1}{c_0} \left\| \theta_2 \Delta t \divergence \testflux \right\|^2
  =: \left| \testflux  \right|_\mathrm{f}^2.
 \end{align*}
 It directly follows that $\mathcal{E}$ is strongly convex wrt.\ $|\cdot|$ with modulus 
$\sigma=1$, and the partial functional derivatives $D_{\displacement}\mathcal{E}$ and 
$D_{\flux}\mathcal{E}$ are uniformly Lipschitz continuous wrt.\ $|\cdot|_\mathrm{m}$ and 
$|\cdot|_\mathrm{f}$ with Lipschitz constants $L_\mathrm{m}=1$ and $L_\mathrm{f}=1$, 
respectively.

 By the H\"older inequality we deduce that 
 \begin{align}
 \label{Eq:proof-aux:beta-1}
  \|\biotcoefficient : \strain(\testdisplacement)\|^2
  &=
  \int_\Omega \left| \biotcoefficient : \strain(\testdisplacement) \right|^2 \, 
\mathrm{d} 
  \boldsymbol x
  \leq 
  \int_\Omega \left| \biotcoefficient : \mathbb{C}^{-1} : \biotcoefficient \right| 
\, \left| \strain(\testdisplacement) : \mathbb{C} : \strain(\testdisplacement) \right| \, 
\mathrm{d} \boldsymbol x\\[1ex]
  \nonumber
  &\leq 
  \left\| \biotcoefficient : \mathbb{C}^{-1} : \biotcoefficient \right\|_\infty
  \llangle \mathbb{C} \strain(\testdisplacement), \strain(\testdisplacement) \rrangle.
 \end{align}
 Hence, it follows that 
 \begin{align*}
  \left|\testdisplacement\right|_\mathrm{m}^2 &\leq \left(1 + 
\frac{\left\|\biotcoefficient : \mathbb{C}^{-1} : \biotcoefficient 
\right\|_\infty}{c_0}\right) \, \left|(\testdisplacement,\testflux)\right|^2.
 \end{align*}
 
 On the other hand, applying the triangle inequality and Young's inequality, and 
balancing the arising constants properly yields that 
 \begin{align*}
  \frac{\theta_1}{c_0} \left\| \theta_2 \Delta t \divergence \testflux \right\|^2
  &\leq 
  \frac{\theta_1}{c_0} \left(1 + \frac{\left\|\biotcoefficient : \mathbb{C}^{-1} : 
\biotcoefficient \right\|_\infty}{c_0}\right) \left\| \theta_2 \Delta t \divergence 
\testflux + \biotcoefficient : \strain(\testdisplacement) \right\|^2\\
  &\hspace{1cm} +
  \frac{\theta_1}{c_0} \left(1 + \frac{c_0}{\left\|\biotcoefficient : \mathbb{C}^{-1} : 
\biotcoefficient \right\|_\infty}\right) \left\| \biotcoefficient : 
\strain(\testdisplacement) \right\|^2.
 \end{align*}
 Together with~\eqref{Eq:proof-aux:beta-1}, we also conclude that
  \begin{align*}
  \left|\testflux\right|_\mathrm{f}^2 &\leq \left(1 + \frac{\left\|\biotcoefficient : 
\mathbb{C}^{-1} : \biotcoefficient \right\|_\infty}{c_0}\right) \, 
\left|(\testdisplacement,\testflux)\right|^2\,.
 \end{align*}

 Thereby, the assumptions of Theorem~\ref{Lemma:AM_convergence} are fulfilled 
and~\eqref{Eq:AM_convergence} is ensured with constants 
$\sigma=L_\mathrm{m}=L_\mathrm{f}=1$ and $\beta_\mathrm{m}=\beta_\mathrm{f} = \left(1 + 
\frac{\left\|\biotcoefficient : \mathbb{C}^{-1} : \biotcoefficient 
\right\|_\infty}{c_0}\right)^{-1}$. Finally, the assertion follows directly, since 
$\mathcal{E}$ is quadratic and $(\displacement^n,\flux^n)$ is a local minimum of 
$\mathcal{E}$ and $|\cdot|$ relates to the second functional derivative of $\mathcal{E}$ 
via~\eqref{Eq:proof-aux:second-derivative-energy}. Therefore, we have that $ 
\mathcal{E}(\displacement^{n,k},\flux^{n,k}) - \mathcal{E}(\displacement^n,\flux^n) = 2 
\left|(\displacement^{n,k} - \displacement^n, \flux^{n,k} - \flux^n)\right|^2$ for all 
$k\geq 0$.
\end{proof}

\begin{rem}[Convergence of $p^{n,k}$]
The convergence of the sequence of pressures $\{p^{n,k}\}_k$ follows now immediately by a 
standard inf-sup argument.
\end{rem}

\begin{rem}[Comparison with quasi-static case]
The final convergence rate in Corollary~\ref{Corollary:Convergence_UD} coincides with the 
one for the undrained split applied to the quasi-static Biot equations for an homogeneous 
and isotropic bulk; cf.~\cite{MW13_BBR}.  In that case, the Biot tensor 
$\biotcoefficient$ 
reduces to $\alpha \mathbf{I}$ for some constant $\alpha\in(0,1]$, and $\mathbb{C}$ is 
defined by the Lam\'e parameters, such that $\biotcoefficient : 
\mathbb{C}^{-1} : \biotcoefficient = \frac{\alpha^2}{K_\mathrm{dr}}$, where 
$K_\mathrm{dr}$ is the drained bulk modulus.
\end{rem}

\section*{Acknowledgement}

This work was supported by the German Academic Exchange Service (DAAD) under the grant ID 
57458510 and by the Research Council of Norway (RCN) under the grant ID 294716.

%
%

\end{document}